\documentclass{amsart}
\usepackage{amsmath, amsthm, amssymb, xypic,graphicx}
\begin{document}
\newcommand{\Z}{\mathbb{Z}}
\newcommand{\Cal}{\mathcal}
\newcommand{\Proj}{\mathbb{P}}
\newcommand{\Fl}{\operatorname{Fl}}
\newcommand{\overbar}{\overline}
\newcommand{\Div}{\operatorname{Div}}
\newcommand{\Nq}{N_{\mathbb{Q}}}
\newcommand{\Mq}{M_{\mathbb{Q}}}
\newcommand{\Q}{{\protect\mathbb{Q}}}
\newcommand{\Vs}{V(\sigma)}
\newcommand{\Nsq}{N(\sigma)_\Q}
\newcommand{\Msq}{M(\sigma)_\Q}
\newcommand{\Qt}{\Q \otimes}
\newcommand{\QDX}{\Q \otimes \Div_T(X)}
\newcommand{\Star}{\operatorname{Star}}
\newcommand{\va}{\widehat{[V(\sigma)]}}
\newcommand{\vz}{[V(\sigma)]}
\newcommand{\di}{\operatorname{div}}
\newcommand{\msqp}{M(\sigma)_\Q^\perp}
\newcommand{\qnsp}{\Qt N_\sigma^\perp}
\newcommand{\ps}{\pi_\sigma}
\newcommand{\pg}{\pi_\gamma}
\newcommand{\pgp}{\pi_{\gamma'}}
\newcommand{\nts}{n_{\tau,\sigma}}
\newcommand{\ntd}{n_{\tau,\delta}}
\newcommand{\nds}{n_{\delta,\sigma}}
\newcommand{\ngs}{n_{\gamma,\sigma}}
\newcommand{\ngps}{n_{\beta,\sigma}}
\newcommand{\ntg}{n_{\tau,\gamma}}
\newcommand{\ntgp}{n_{\tau,\beta}}
\newcommand{\g}{\gamma}
\newcommand{\gp}{{\gamma'}}

\newcommand{\s}{\sigma}
\newcommand{\nttg}{{\tilde n}_{\tau,\gamma}}
\newcommand{\nttgp}{{\tilde n}_{\tau,\beta}}
\newcommand{\la}{\langle}
\newcommand{\ra}{\rangle}
\newcommand{\Nst}{N(\s)_\tau}
\newcommand{\Mst}{M(\s)_\tau}
\newcommand{\Ngpt}{N(\beta)_\tau}
\newcommand{\Mgt}{M(\g)_\tau}
\newcommand{\zq}{Z_*(X)_\Q}
\newcommand{\zk}{Z_*(X)_\K}
\newcommand{\ntr}{n_{\tau_j,\rho_j}}
\newcommand{\nttr}{\tilde n_{\tau_j,\rho_j}}
\newcommand{\nttrk}{\tilde n_{\tau_k,\rho_k}}
\newcommand{\qyn}{\Q [Y_1,\dots,Y_n]}
\newcommand{\qyr}{\Q [Y_1,\dots,Y_n]}
\newcommand{\kyn}{\K [Y_1,\dots,Y_n]}
\newcommand{\kyr}{\K [Y_1,\dots,Y_n]}
\newcommand{\Nr}{N_\mathbb{R}}

\newcommand{\K}{\mathbb{K}}
\newcommand{\Kt}{\K \otimes}
\newcommand{\avol}{\operatorname{vol}}

\newcommand{\p}{\pi}
\newcommand{\product}{\prod}
\newcommand{\comp}{\circ}

\newcommand{\F}{\mathbb{F}}
\newcommand{\psp}{\pi_{\s'}}
\newcommand{\Ps}{\Psi}
\newcommand{\om}{\omega}

\newcommand{\Os}{\Cal{O}}
\newcommand{\Td}{\operatorname{Td}}
\newcommand{\ch}{\operatorname{ch}}
\newcommand{\ls}{\mathcal{L}}
\newcommand{\Gr}{\operatorname{Gr}}
\newcommand{\Hom}{\operatorname{Hom}}

\newcommand{\bP}{{\underline{ \Psi}}}

\newtheorem{theorem}{Theorem}[section]
\newtheorem{lemma}{Lemma}[section]
\newtheorem{prop}{Proposition}[section]
\newtheorem{cor}{Corollary}
\newtheorem*{cor*}{Corollary}
\theoremstyle{remark}
\newtheorem*{rem}{Remark}
\newtheorem*{con}{Convention}
\theoremstyle{definition}
\newtheorem*{ex}{Example}
\title{Cycle-Level Intersection Theory for Toric Varieties}

\author{Hugh Thomas}
\address{Fields Institute\\ 222 College St.\\
Toronto, ON \\M5T 3J1}
\email{hthomas@fields.utoronto.ca}


\begin{abstract} This paper addresses the problem of constructing a 
cycle-level intersection theory for toric varieties.  
We show that by making one global choice, 
we can determine a cycle representative
for the intersection of an equivariant Cartier divisor with an invariant 
cycle on a toric variety.  For a  toric variety
defined by a fan in $N$, the choice consists of giving an 
inner product or a complete flag for $M_\Q=
\Qt \Hom(N,\mathbb{Z})$, or more
generally giving for each cone $\s$ in the fan a linear subspace of 
$M_\Q$ complementary to $\s^\perp$, satisfying certain compatibility
conditions.  
We show that these intersection cycles have properties analogous to the 
usual intersections modulo rational equivalence.  
If $X$ is simplicial (for instance, if $X$ is non-singular), 
we obtain a commutative ring structure
to the invariant cycles of $X$ with rational 
coefficients.  This ring structure determines cycles representing 
certain characteristic classes of the toric variety.  
We also discuss
how to define intersection cycles that require no choices, 
at the expense of increasing
the size of the coefficient field.

\end{abstract}

\maketitle

\section{Introduction}

The basic building block of intersection theory is the intersection of 
a Cartier divisor and a cycle, which is defined to be a cycle modulo 
rational equivalence.  In this paper, we restrict our attention to
toric varieties, and
we ask whether there is a way to define the intersection of an equivariant
Cartier divisor and an invariant cycle as a cycle, rather than a
cycle modulo rational equivalence.  This cannot be done naturally; a choice
must be made.  The goal of this paper is to show that a simple 
initial choice allows us to define an action of equivariant Cartier divisors
on invariant cycles in a consistent manner.  From this, a cycle-level
intersection theory can be developed.  For simplicial toric varieties, this
gives rise to an intersection product on the level of cycles, analogous
to the usual intersection product on cycles modulo rational equivalence.

First, we describe the properties which our action posesses.  Let $X$ be a 
toric variety, $D$ and $E$ equivariant Cartier divisors on $X$, $z$ an 
invariant 
cycle of $X$ with coefficients in $\Q$, and $\overbar z$ the image of 
$z$ in the Chow group of cycles modulo rational equivalence.  Let $D\cdot{}$
denote either our action of equivariant Cartier divisors on cycles, or the usual action
of Cartier divisors on cycles modulo rational equivalence, depending on context.
Then the action we define has the following properties:
\begin{enumerate}
\item $\overbar {D\cdot z}=D\cdot \overbar z$, i.e., when we pass to rational
equivalence, our intersection cycle agrees with the usual intersection defined
on cycles modulo rational equivalence.
\item If $D$ and $z$ intersect properly, then there is a naturally well-defined
cycle-level intersection, namely $[D|_z]$, the cycle corresponding to the
restriction of $D$ to $z$.  In this case, $D\cdot z=[D|_z]$.  
\item $D \cdot (E \cdot z)=E\cdot (D \cdot z)$.  
\end{enumerate}

Now, we briefly describe the choice that must be made in order to define
our action.  Let $N$ be a free abelian group of rank $n$.  Let $\Sigma$ be
a fan in $N_\Q=N\otimes \Q$, that is, a collection of strongly convex 
rational polyhedral cones in $N_\Q$, such that any face of a cone $\sigma$
in $\Sigma$ is again in $\Sigma$ and 
 the intersection of any two
cones in $\Sigma$ is a face of both and is again in $\Sigma$.  (For 
further details about toric varieties, see Section 2 and 
\cite{F2}.)
Let $M=\Hom(N,\Z)$.  Each $\sigma \in \Sigma$ defines a subspace
$\sigma^\perp$ of $M_\Q=\Qt M$.  To define our action, we need a function $\Psi$
associating to each cone $\sigma$ of $\Sigma$ a subspace $\Psi(\sigma)$
with the following two properties:
\begin{enumerate}
\item For each $\sigma \in \Sigma$, $\Psi(\sigma)$ and $\sigma^\perp$ are
complementary subspaces of $M_\Q$.
\item If $\sigma \in \Sigma$ and $\tau$ is a face of $\sigma$, then
$\Psi(\sigma)$ contains $\Psi(\tau)$.  
\end{enumerate}

There are two easy ways to define a map $\Psi$ having these properties, 
starting with either an inner product on $M_\Q$ or a complete flag in $M_\Q$
which is generic in a suitable sense with respect to $\Sigma$.  

Sections 2--6 are concerned with constructing our action and showing that it
has the properties we have already described.  In section 7, we consider the
situation where we have two toric varieties and a proper equivariant map between
them $f:X'\rightarrow X$.  In this case, if we choose $\Psi$ for $X$ and 
$\Psi'$ for $X'$ in compatible ways, we have an analogue of the projection
formula from intersection theory: if $z$ is an invariant cycle
on $X'$ and $D$ is an equivariant Cartier divisor on $X$,  then 
$$f_*(f^*(D)\cdot z)=D\cdot f_*(z).$$
In sections 8--12 we consider the case where $\Sigma$ is simplicial.  In this
case we obtain a commutative ring structure on 
the invariant cycles of $X$ with rational coefficients.  In sections 13--15 we discuss techniques for
performing computations using our action.  Section 16 concerns how we
can use the actions already constructed
 to produce actions of equivariant Cartier divisors
on invariant cycles independant of any choice, at the expense of enlarging 
the coefficient field.  Section 17 considers how one can go about 
constructing cycles which represent characteristic classes of a 
toric variety.

\section{Preliminaries}

We begin by recalling necessary background about toric varieties.  For further
explanation and proofs, see \cite{F2}.  

Let $N$ be a free abelian group of rank $n$, and $M=\Hom(N,\Z)$ its dual.  
We will
denote $\Qt N$ by $\Nq$.  In general, a subscript $\Q$ will mean ``tensor
by $\Q$'', and similarly for other fields. 
  
Let $\Sigma$ be a fan in $N$, and $X=X(\Sigma)$ the associated
toric variety.  (To define $X(\Sigma)$ we have to choose a ground field, but
this choice is of no importance for our purposes.)
If $\tau$ and $\sigma$ are two cones of $\Sigma$, 
$\sigma \succ \tau$ means that $\tau$ is a proper face of $\sigma$, and $\sigma
\rightarrow \tau$ means that $\tau$ is a maximal proper
face of $\sigma$, also known as a facet.  We shall
say that $\Sigma$ is affine if it consists of a single cone and its
faces.  

For $\sigma \in \Sigma$, $N_\sigma$
is the subgroup generated by $\sigma \cap N$, and  $N(\sigma) = N/N_\sigma$.  
$M(\sigma)$ is the dual of $N(\sigma)$, so naturally a sublattice of $M$;
in fact, $M(\sigma)=M \cap \sigma^\perp$.  
$M_\sigma = M/M(\s)$ is dual to $N_\sigma$.    
For $\sigma \prec \tau$, we set $M(\sigma)_\tau = M(\s)/M(\tau)$,
and $N(\sigma)_\tau = N_\tau/N_\sigma$.  These two are also dual.  
All the duality pairings mentioned in this paragraph are denoted $\langle 
\cdot,\cdot \rangle$.  

Let $D \in \Div_T(X)$, the equivariant Cartier divisors on $X$ 
(from now on, we 
shall omit the word equivariant).  
Define $m_\sigma \in M_\sigma$ by requiring that $\chi^{m_\sigma}$
be a local equation for $D$ on $U_\sigma$.  (Frequently in the study of 
toric varieties it is preferred to take the local equations for $-D$.  To 
follow that convention would introduce a number of negative signs, and 
so we have
chosen not to.)  By an abuse of terminology, we will refer to the 
collection of $m_\sigma$ associated to a divisor $D$ as its local equations.
The $m_\sigma$ satisfy certain agreement conditions: if $\tau \prec \sigma$,
the image of $m_\sigma$ under the natural map from $M_\sigma$ to $M_\tau$
is $m_\tau$.  Conversely, given a collection of $m_\s$ satisfying these 
conditions, they determine a Cartier divisor.  

We shall also want to consider $\Q$-Cartier divisors, by which we mean 
rational multiples of Cartier divisors.  More precisely, we will think
of a $\Q$-Cartier divisor as a collection of $m_\s \in \Qt M_\sigma$,
satisfying the same agreement conditions as above for Cartier divisors.

Given $\sigma \in \Sigma$ a cone of dimension $n-k$, there is an associated 
$k$-dimensional subvariety of $X$, which is denoted $V(\sigma)$.  
$V(\sigma)$ is also a toric variety, and is given by a fan in $N(\sigma)$.
Its fan is 
$\Star(\sigma)=\{ (\tau+N_\sigma) /N_\sigma| \tau \in \Sigma, 
\tau \succeq \sigma\}$. 

The group of invariant $k$-cycles, which we will denote $Z_k(X)$, 
is the free abelian group on the set of $[V(\sigma)]$, where $\sigma$ 
ranges over the cones of dimension $n-k$ in $\Sigma$. From now on, we omit
the word ``invariant'': all our cycles are understood to be torus invariant.
Let $\rho_1 , \dots , \rho_r$
be the rays of $\Sigma$.  
Denote $V(\rho_i)$ by $D_i$.  

Denote by $A_k(X)$ the $k$-th Chow group of $X$.  It 
has a presentation as the quotient
of $Z_k(X)$ by the invariant $k$-cycles which are rationally 
equivalent to zero.  (See \cite{FS} for more details.)

There is a map from Cartier divisors to $Z_{n-1}(X)$.  If $D$
is a Cartier divisor,  the associated element of
$Z_{n-1}(X)$ is $[D]=\sum_{i=1}^r 
\la m_{\rho_i},v_i\ra D_i$, where $v_i$ is the first lattice point of 
$N$ along $\rho_i$ and $m_{\rho_i}$ is the local equation for $D$ on $\rho_i$.
The same 
definition gives a map from $\Q$-Cartier divisors to $Z_{n-1}(X)_\Q$
which extends the previous map.

We shall also be interested in this map from $\Q$-Cartier divisors to 
invariant cycles when the divisors are on $V(\sigma)$ rather than
$X$.  Let $D$ be a $\Q$-Cartier divisor on $V(\sigma)$ and let $\sigma$ be 
a maximal proper 
face of $\tau$.  Then the 
coefficient of $[V(\tau)]$ in $[D]$ is $\la m_\tau,\nts \ra$, where $\nts$ is the
generator of the (1-dimensional) image of $\tau$ in $N(\sigma)$, i.e.,~the
first lattice point of $N(\sigma)$ along the ray in $N(\sigma)$ corresponding
to $\tau$.  

\section {A Little Linear Algebra}

Let $<$ denote the relation of being a vector
subspace.  We establish two lemmas.

\begin{lemma} Let $V\geq A \geq B$. Let $A'$ and $B'$ be complementary 
subspaces to $A$ and $B$ in $V$  such that $B' \geq A'$.  Then 
\begin{enumerate}
\item[(i)]  There is
a canonical isomorphism between $A/B$ and $B'/A'$.

\item[(ii)] There is a canonical map from $V/B$ to $A/B$.
\end{enumerate}
 \end{lemma}

\begin{proof} This follows from the fact that
$V\cong B \oplus (A \cap B') \oplus A'$, $A \cong B \oplus (A \cap B')$.  
\end{proof}

\begin{lemma} Let $V \geq A \geq B \geq C$.  Let $A'$, $B'$, and $C'$ be 
complementary 
subspaces to $A$, $B$, and $C$, such that $V \geq C' \geq B' \geq A'$.  
Then the following
diagram commutes, where the horizontal arrows are projection, and the vertical
arrows are the maps from Lemma 3.1.
$$
\diagram
V/C \rto \dto & V/B \dto \\
A/C \rto & A/B
\enddiagram
$$
\end{lemma}

\begin{proof} The proof is similar to that of the previous lemma.  Write
  $V \cong C \oplus (B\cap C') \oplus (A \cap B') \oplus A'$, and 
similarly for $A$ and $B$.    
\end{proof}

\section {Definition of the Action and Basic Properties}

\begin{con} I shall frequently need to refer to sub- and quotient 
groups of $M$ tensored by $\Q$, and very seldom to the groups themselves.  
Therefore, I shall drop the $\Q \otimes{}$ from my notation: whenever I refer
to $M$, $M(\s)$, $M_\tau$, $M(\s)_\tau$, etc., the $\Q \otimes{}$ is to be 
understood.  Note that this does not apply to sub- and quotient groups of
$N$.  \end{con}

To define the action, we fix a map $\Psi$ taking cones
of $\Sigma$ to linear subspaces of $ M$ of the same dimension,
satisfying the following properties
\begin{enumerate}
\item[(i)] $\Psi(\s)$ and $M(\s)$ span $M$ (or equivalently, 
$\Psi(\s) \cap (M(\s))= \{0\}$).

\item[(ii)] If $\s \succ \tau$, $\Psi(\s) > \Psi(\tau)$.  
\end{enumerate}
If $\Psi$ satisfies these two properties, we  call $\Psi$ a choice of complements for $\Sigma$.

We immediately give possible constructions of $\Psi$.  

\begin{ex}[Inner Product] Fix an inner product on 
$ M$.  Define $\Psi(\s)$ to be the orthogonal complement of $\s^\perp$
in $ M$.  \end{ex}

\begin{ex}[Generic Flag] Fix a complete flag in $M$,
$0=F_0 < F_1<\dots F_n=M$, which is generic with respect to $\Sigma$ 
in the sense that if 
$\sigma$ is $k$-dimensional, then $\s^\perp\cap F_k=\{0\}$. Define
$\Psi(\s)=F_k$ for all $\s$ of dimension $k$. \end{ex}


Let $D$ be a $\Q$-Cartier divisor and $\sigma$ a cone of dimension $n-k$ 
in $\Sigma$.  
We proceed to define the action of
$D$ on $[V(\s)]$.

The decomposition $M= M(\s) \oplus \Psi(\s)$ determines a projection
$\ps: M \rightarrow M(\s)$.
Also, for any $\tau \succ \s$, writing $M_\tau = M/M(\tau)$, $M(\s)_\tau=
M(\s)/M(\tau)$, we can apply Lemma 3.1 to get a map, denoted the same way, 
$\ps:M_\tau \rightarrow M(\s)_\tau$.  
(Observe that in the case of the inner product action,
the map $\ps$ is orthogonal projection from $M$ to $M(\s)$.)

For $\tau \in \Sigma$, let the local equation of 
$D$ on $\tau$ be $m_\tau$.  For $\tau \succ \s$, 
define 
$\overbar m_{\overbar\tau} \in M(\s)_\tau$ by $\overbar m_{\overbar\tau} = 
\ps(m_\tau)$.
We verify that these $\overbar m_{\overbar\tau}$ form local equations for a 
$\Q$-Cartier divisor on the fan of $V(\s)$ in $M(\s)$.  
Suppose $\tau \succ \gamma \succ \sigma$.
We need to show that $\overbar m_{\overbar\tau}$ maps to $\overbar 
m_{\overbar\g}$ under the 
natural map from $M(\s)_\tau$ to $M(\s)_\g$.  But this follows
immediately from Lemma 3.2 and the fact that $m_\tau$ goes to $m_\g$ under the map from
$M_\tau$ to $M_\g$. Thus, the $\overbar m_{\overbar \tau}$ form local 
equations for a
$\Q$-Cartier divisor, which we denote $D_\s$.

So, define  
$D \cdot [V(\sigma)]= [D_\sigma]$.  
Explicitly, this definition is:
$$D \cdot [V(\s)] = \sum_{\tau \rightarrow \s} \la \ps(m_\tau),n_{\tau,\s} \ra
[V(\tau)].$$
For any $\Q$-Cartier divisor and any $k$, this gives us a map
$Z_{k}(X)_\Q \rightarrow Z_{k-1}(X)_\Q$, as desired.

\begin{ex}[$\mathbb{C}^n$ with the standard inner product] 
Let $e_1 , \dots, e_n$ be a basis of $N$,
and $e_1^*,\dots,e_n^*$ the dual basis of $M$. Put the usual inner product
on $M$ so that the $e^*_i$ form an orthonormal basis..  
Let us consider the example where $\Sigma$
is the fan consisting of the positive orthant $\s$ and its faces, so 
$X(\Sigma) \cong \mathbb{C}^n$.  
Let $\rho_i$ be the ray in the positive $e_i$ direction, and $D_i=
V(\rho_i)$.  Then $D_i \cdot V(\tau)=0$ if $\rho_i \prec \tau$, and 
$D_i \cdot V(\tau)=V(\tau+\rho_i)$ otherwise.  
\end{ex}

Given $D$ a $\Q$-Cartier divisor, we say the $D$ intersects $V(\s)$ properly 
if when $D$ is written as a sum of codimension one irreducible subvarieties,
no subvarieties containing $[V(\s)]$ occur with non-zero coefficients.
Equivalently, $D$ intersects $[V(\s)]$ properly if it restricts to a 
$\Q$-Cartier
divisor on $V(\s)$.  If $m_\s$ is
the local equation of $D$ on $\s$, $D$ intersects $V(\s)$ properly iff 
$m_\s=0$.

We now examine $D \cdot V(\sigma)$ when $D$ intersects $V(\sigma)$ 
properly. 

\begin{prop} Let $D$ be a $\Q$-Cartier divisor, and suppose it 
intersects $[V(\s)]$ properly.  Then $D_\s$ is the restriction of $D$ to 
$V(\s)$, so $D \cdot [V(\sigma)]$ is the cycle 
corresponding to the restriction of $D$ to $V(\s)$.  \end{prop} 
\begin{proof}
Since $m_\s=0$, the image of $m_\tau$ in $M_\s$ is $0$, so
$m_\tau \in M(\s)_\tau$, and thus $\ps(m_\tau)=m_\tau$.  Thus, the local
equations for $D_\sigma$ are just the local equations for $D$.
Thus we see that $D_\sigma$ is the restriction of $D$ to 
$V(\sigma)$, and from the definition of the action it follows that
$D\cdot [V(\sigma)]$ is just the cycle associated to the restriction
of $D$ to $V(\sigma)$. \end{proof}

If $X$ is a general algebraic variety, and $D$ a Cartier divisor on $X$,
$D$ induces a map from $A_{k}(X)$ to $A_{k-1}(X)$, also denoted $D \cdot{}$.  
(See \cite{F1} for details.)  If $z \in Z_k(X)_\Q$, denote by $\overline{ z}$ its image
in $A_k(X)_\Q$.  
We now show that our action agrees with the usual action of Cartier divisors
on Chow groups once we pass to rational equivalence.

\begin{prop}
If $z \in Z_k(X)_\Q$, the image of $D \cdot z$ in $A_{k-1}(X)_\Q$ equals
$D \cdot \overline{ z}$.  

\end{prop}

\begin{proof}
$D\cdot \overline{[V(\s)]}$ is defined by translating $D$ 
by a principal divisor so that it intersects $V(\sigma)$
properly, then taking the class modulo rational equivalence of the
intersection.  Observe that $D_\sigma$ is the restriction to $V(\s)$ of
a translate of $D$ by a principal divisor, so 
$\overline{[D_\s]}=D\cdot\overline{[V(\s)]}$, as 
desired.  
\end{proof}

\section{All Computations Can be Done on Affine Varieties}

Recall that $X$ is covered by affine open sets $U_\s$ with $\sigma
\in \Sigma$, where $U_\s$ is the toric variety corresponding to the fan 
consisting of $\sigma$ and all its faces.   This covering is useful because it allows us in some
cases to localize the questions we are interested in to $U_\s$ and deal with
them there, where we don't need to worry about the global structure of the fan.

We recall some useful facts about the covering of $X$ by
the $U_\s$.  If $\tau \prec \s$, then we can perform the 
construction of $V(\tau)$ in the toric variety $U_\sigma$.  Let us denote this 
subvariety of $U_\s$ by $V(\tau)_\sigma$.  Then $V(\tau)_\sigma = V(\tau) \cap
U_\sigma$  is a dense
open subset of $V(\tau)$.  (See \cite{F2} for more details.)

\begin{lemma} Let $\tau \rightarrow \s$.  
Let $D$ be a $\Q$-Cartier divisor on $X$.  
Then the coefficient of $[V(\tau)]$
in $D \cdot [V(\s)]$ equals the coefficient of $[V(\tau)_\tau]$ in
$D |_{U_\tau} \cdot [V(\s)_\tau]$.  \end{lemma}

\begin{proof} 
Let $m_\tau$ be the local equation of $D$ on $\tau$.  
The local equation of the restriction of $D$ to $U_\tau$ is also $m_\tau$.  
Thus,
the  coefficients in question both equal  
$\la\ps(m_\tau), \nts \ra$.
\end{proof}

The analogous result for the effect of a composition of the action of 
multiple Cartier divisors is also true.  For simplicity, if $D$ and $E$
are $\Q$-Cartier divisors and $z$ is a cycle, 
instead of
writing $D \cdot (E \cdot z)$, we write $D \cdot E \cdot z$. 

\begin{cor*} Let $\tau \succ \s$.  Let $E_1,\dots,E_s$ be 
$\Q$-Cartier divisors on $X$.  Then 
the coefficient of 
$[V(\tau)]$ in $E_1 \cdot \cdots  E_s \cdot [V(\s)]$ equals the 
coefficient of $[V(\tau)_\tau]$ in $E_1 |_{U_\tau} \cdot \cdots  
E_s |_{U_\tau} \cdot
[V(\s)_\tau]$.  \end{cor*}

\begin{proof} This follows by a repeated application of Lemma 5.1.  
\end{proof}

In fact, computing the coefficient of a cycle $[V(\tau)]$ as above
can be reduced to computing the coefficient of a cycle $[V(\tau')]$ 
where $\tau'$ is of maximal dimension in a lattice $N'$.  

Let $\tau$ be a cone of $\Sigma$.  
Let $N'=N_\tau$, and let $\Sigma'$ be the fan in $N'$ consisting of $\tau$ and 
its faces.  For $\g \preceq \tau$, let $\g'$ denote the corresponding cone of 
$\Sigma'$.  
Let $E_1,\dots,E_s$ be $\Q$-Cartier divisors on $X$, and 
for $1\leq i \leq s$ let
$E_i'$ be the $\Q$-Cartier divisor on $X(\Sigma')$ having the same 
local equations on $\g'$ as $E_i$ does on $\g$, for all $\g \preceq \tau$.  
Then we have the following result:

\begin{lemma} The coefficient of $[V(\tau)]$ in $E_1 \cdot \cdots  E_s
\cdot [V(\s)]$ equals the coefficient of $[V(\tau')]$ in 
$E_1' \cdot \cdots  E_s'
\cdot [V(\s')]$. \end{lemma}

\begin{proof} As in the proof of Lemma 5.1, the explict expression for the 
two coefficients is the same.  \end{proof}

\section{Commutativity}
{
\renewcommand{\gp}{\beta}
\renewcommand{\pgp}{\pi_{\beta}}

This section is devoted to proving the following commutativity result:

\begin{theorem} Let $D$ and $E$ be $\Q$-Cartier divisors on $X$,
and $\sigma$ a cone of $\Sigma$ of dimension $n-k$.  Then $D\cdot E\cdot
[V(\sigma)] = E \cdot D \cdot [V(\sigma)]$. \end{theorem}

\begin{proof}
$D \cdot E \cdot [V(\sigma)]$
is a linear combination of  $[V(\tau)]$ where $\tau$ ranges over cones of
dimension $n-k+2$.  More explicitly, if the local equations for $E$ are
$e_\delta$,
$$E\cdot V(\sigma) = \sum_{\delta \rightarrow \sigma} 
\la\ps(e_\delta), \nds  \ra [V(\delta)]$$
and, if the local equations for $D$ are $d_\tau$,
$$D \cdot E \cdot [V(\sigma)]=
\sum_{\tau \rightarrow \delta}\sum_{\delta 
\rightarrow \sigma}\la \pi_\delta(d_\tau),\ntd\ra\la \pi_\s (e_\delta),\nds
\ra [V(\tau)].
$$
Suppose we wish to compute the coefficient of some particular $[V(\tau)]$
in the above expression. 
From convex geometry, we know that there are exactly two cones $\delta$
satisfying $\tau \rightarrow \delta \rightarrow \sigma$.  Call these cones
$\gamma$ and $\beta$.  Then the coefficient of $[V(\tau)]$ in $D\cdot E\cdot 
[V(\s)]$ is  
\begin{equation} \la \pi_\g(d_\tau),\ntg\ra\la \pi_\s (e_\g),\ngs
\ra + \la \pi_\gp(d_\tau),\ntgp \ra\la \pi_\s (e_\gp),\ngps
\ra. 
\end{equation}

To prove that $D\cdot E\cdot [V(\s)] = E\cdot D\cdot 
[V(\s)]$, we need to check that the coefficient of $[V(\tau)]$ 
in the two expressions
is the
same for all $\tau$.  
Thus, to show commutativity, it suffices to show that the 
expression (1) is the same after we change all the $d$'s to $e$'s and vice
versa.

Our first goal is to modify (1) so that all the pairings in it are between
the same pair of
spaces, $M(\s)_\tau$ and $\Qt N(\s)_\tau$.  By replacing $\pi_\s(e_\g)$ and
$\pi_\s(e_\gp)$ by $\pi_\s(e_\tau)$, which we may do by the 
compatibility conditions for local equations, we have put the second and 
fourth pairings of (1) into the desired form.  

Figure 1 shows an example where $\Psi$ is induced from the standard inner
product. 

\begin{figure}[htb]
\includegraphics*{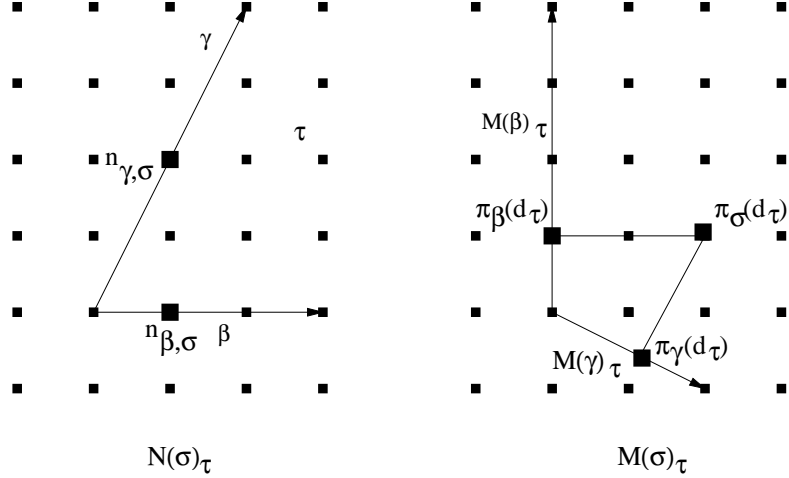}
\caption{An illustration of $N(\s)_\tau$ and $M(\s)_\tau$}
\end{figure}

We move on to consider the first and third pairings.  
Both $\ntg$ and $\ntgp$ live in
quotients of $N(\s)_\tau$, so we choose liftings of them to 
$\Qt N(\s)_\tau$: define $\nttg \in 
\Qt N(\s)_\tau$ to be the lifting of $\ntg$ which lies in $\gp$.  Similarly,
define $\nttgp$ to be the lifting of $\ntgp$ which lies in $\g$.   
(See Figure 2.)
Thus we can write the coefficient of $[V(\tau)]$ in $D\cdot E \cdot [V(\s)]$
as 
$$ \langle\ps(e_\tau), \ngs \rangle \langle \pg(d_\tau),\nttg \rangle +
\langle\ps(e_\tau),\ngps \rangle\langle \pgp(d_\tau),\nttgp \rangle $$
where all the pairings are between $M(\s)_\tau$ and $\Qt N(\s)_\tau$.

We would now like to analyze $\nttg$ and $\nttgp$ further.  
It is clear by construction that there are constants 
$c$ and $c'$ such that $c\nttg =\ngps$ and $c' \nttgp=\ngs$.  We need
to show that $c=c'$.  We do this by the following lemma:

\begin{lemma} The following equalities hold:
\begin{eqnarray*} & & {[N(\s)_\tau:L]}\nttg =  \ngps \\
& &  {[N(\s)_\tau:L]}\nttgp= \ngs. \end{eqnarray*}
\end{lemma}

\begin{proof}
We prove the first equality.

Let $W$ denote the fundamental domain of $L$ whose vertices are
$0$, $\ngs$, $\ngps$, and $\ngs + \ngps$, and we consider $W$ not to include
any points from its left or top sides; in particular 
the unique point of $L$ contained in $W$ is $\ngps$.  (See Figure 2.)
The number of 
points of $\Nst$ in $W$ is then $[\Nst:L]$.  I will proceed to show that 
the number of points is also $c$.  
\begin{figure}[htb]
\begin{center} \includegraphics*{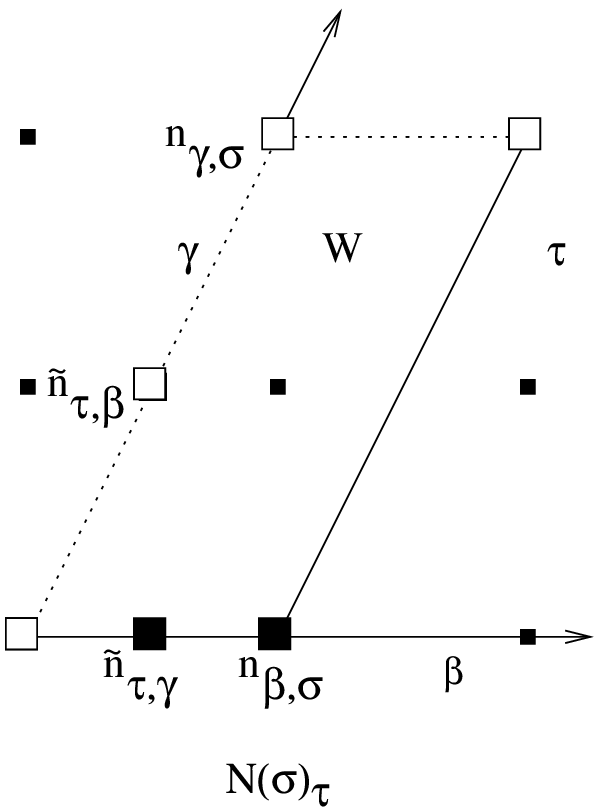} \end{center}
\caption{An illustration of $N(\s)_\tau$ showing $W$}
\end{figure}

Consider the projection $\xi:\Nst \rightarrow N(\g)_\tau = \Nst/\mathbb{Z}\cdot \ngs
=\Z\cdot n_{\tau,\g}$. 
$N(\g)_\tau$ is a one-dimensional lattice,
so we can put a linear order on it by saying that $n_{\tau,\gamma}$ designates
the positive direction.  Now
$\xi$ defines a 
bijection between $W \cap \Nst$ and 
$\{ x \in N(\g)_\tau | 0 <x \leq \xi(\ngps) \}$.
The size of the first set is $[\Nst:L]$, while
the size of the second set is $\xi(\ngps)/\ntg = c$.  This completes the proof
of the lemma. \end{proof}

Thus the coefficient of
 $[V(\tau)]$ in $D\cdot E \cdot [V(\s)]$
is 
$$ \!\!\!\!\!\frac{1}{[\Nst:L]} \bigl(
 \langle \ps(e_\tau),\ngs \rangle \langle \pg(d_\tau),\ngps \rangle +
\langle\ps(e_\tau),\ngps \rangle\langle \pgp(d_\tau),\ngs \rangle 
\bigr).$$

Coordinatize $M(\s)_\tau$ by $\la \cdot,\ngs \ra$, $\la \cdot,\ngps \ra$.
Let $\ps(d_\tau)=(d_1,d_2)$, $\ps(e_\tau)=(e_1,e_2)$.  There are constants
$a$ and $b$ such that
$\pg((x,y))=(0,ax+y)$, while $\pgp((x,y))=(x+by,0)$.  Then the expression above 
simplifies to: 
$$\frac{1}{[\Nst:L]} (ad_1e_1 +d_1e_2 +d_2e_1 +bd_2e_2)$$
which is unchanged if $d$'s and $e$'s are interchanged, as desired.
\end{proof}
}

\section{Proper Maps and Push Forwards}

Let $N'$ be a free abelian group of rank $n'$, and $N$ a free abelian
group of rank $n$.  Let $\Sigma'$ and $\Sigma$ be fans in $N'$ and $N$ 
respectively.  The set of fans in lattices is made into a category by
defining the morphisms from $(N',\Sigma')$ to $(N, \Sigma)$ 
to be the  
group homomorphisms from $N'$ to $N$ such that the image of every cone in 
$\Sigma'$
is contained in some cone in $\Sigma$.  Such a map of fans induces an 
equivariant morphism of toric varieties.  
From now on, whenever we refer to a map
between toric varieties, we mean 
one which is induced from a map of fans.  

A proper map of toric varieties has a characterization in terms of the
map of fans: given a map $\phi$ from $(N', \Sigma')$ to $(N,\Sigma)$, 
the corresponding map of toric varieties is proper iff $\phi^{-1}(|\Sigma|)
=|\Sigma'|$.  If $\phi$ induces a proper map of toric varieties, we say
that $\phi$ is a proper map.  

In general, given a proper map $f: X' \rightarrow X$ of algebraic varieties, 
there is a 
push forward map $f_*$ on cycles. If $Z$ is a $k$-dimensional subvariety of
$X'$, let $W$ be the closure of the image of $Z$.  If the dimension of
$W$ is less than that of $Z$, $f_*([Z])=0$.  If the dimensions are the same,
$f_*([Z])=[k(Z):k(W)][W]$, where $k(\cdot)$ denotes the function field.
(See \cite{F1} for details.)

In the case of toric varieties and invariant subvarieties, this has a 
 translation to lattices.  Again, suppose $f$ is a proper map from
$X'$ to $X$, toric varieties given by fans $\Sigma'$ in $N'$ and $\Sigma$
in $N$, and $\phi$
is the corresponding map from $(N',\Sigma')$ to $(N,\Sigma)$.  Suppose
for convenience that $\Sigma$ is full-dimensional in $N$, that is, that
$\Sigma$ is not contained in any proper subspace of $N$.  
Let $\sigma'$ be a cone in $\Sigma'$.  By the defininition of a map, 
$\phi(\sigma')$ is
contained in some cone of $\Sigma$.  Taking the intersection of all
cones of $\Sigma$ containing $\phi(\sigma')$, 
we find a minimal cone $\sigma \in \Sigma$ containing $\phi(\sigma')$.
Then if the codimension of $\sigma$ is less than the codimension of $\s'$,
$f_*([V(\s')])=0$.  If the codimensions are equal, $f_*([V(\sigma')])=
[M'(\sigma'):\phi^*(M(\sigma))][V(\sigma)]$, where $\phi^*$ denotes the
induced map from $M$ to $M'$, and $[ \cdot : \cdot]$ denotes the
index of the second lattice in the first.  (Here we mean
the actual lattices, without tensoring by $\Q$.)

In the general context of a proper map $f:
 X' \rightarrow X$ of algebraic varieties, there are two important facts
about $f_*$.  Firstly, it respects rational equivalence, so it passes 
to Chow groups, and secondly, given $\overline z$ an element of a Chow group of 
$X'$ and
a Cartier divisor $D$ on $X$, $f_*(f^*(D) \cdot \overline z)=D \cdot f_*(
\overline z)$.  

We proceed to prove an analogous identity for our action of $\Q$-Cartier divisors
on invariant cycles.  

\begin{theorem} Let $f:X' \rightarrow X$ be a proper map of
toric varieties, with corresponding fans $\Sigma'$ in $N'$ and $\Sigma$ in 
$N$ and corresponding map $\phi: N' \rightarrow N$.  
  Let $\Ps'$ and $\Ps$ be choices of complements for 
$\Sigma'$ and $\Sigma$, such that 
if $\phi(\s') \subset \s$ with 
$\s'$ and $\s$ of the same codimension, then $\phi^*(\Ps(\s)) \subset 
\Ps (\s')$.
Let $D$ be a $\Q$-Cartier divisor on X.  Let 
$\sigma' \in \Sigma'$.  Then
$$ f_*(f^*(D)\cdot [V(\s')]) = D \cdot f_*([V(\s')]). 
$$
\end{theorem}

\begin{proof}
For $\g' \in \Sigma'$ let
$c(\g')$ denote the smallest cone in $\Sigma$ containing $\g'$.  Let
$P$ denote the set of $\g'\in\Sigma'$ such that the codimension of 
$\g'$ equals the codimension of $c(\g')$.  (Note that $\g' \in P$ iff
$f_*([V(\g')]) \ne 0$.)

Let the dimension of $\sigma'$ be $n'-k$.  
The following lemma describes the possible arrangements of $\s'$, 
the cones over $\s'$, and smallest cones in $\Sigma$ containing them.  

\begin{lemma}  Let $\Sigma$ in $N$ and $\Sigma'$ in $N'$ be 
two fans, and let $\phi$ be a proper map from $N'$ to $N$.  
Let $\s'$ be a codimension-$k$ cone of $\Sigma'$.  Then
one of the three following situations must occur:
\begin{enumerate}
\item[(i)] $\s' \not \in P$, 
and for all $\tau' \rightarrow \s'$, $\tau' \not \in P$.

\item[(ii)] $\s' \not \in P$, there are two cones $\tau', \g' \in P$ such that
$\tau' \rightarrow \s'$, $\g' \rightarrow \s'$, $c(\s')=c(\tau')=c(\g')$, and
for any other $\delta' \rightarrow \s'$, $\delta' \not \in P$.  

\item[(iii)] $\s' \in P$. 
For each $\tau \rightarrow c(\s')$, there is exactly
one $\tau' \rightarrow \s'$ such that $c(\tau')=\tau$.  Any other $\g' 
\rightarrow \s'$ ( i.e.~one which doesn't correspond to some $\tau
\rightarrow c(\s)$) is not contained
in $P$. 
\end{enumerate} \end{lemma}

The three situations are shown in Figure 3.  In order to represent 
a three-dimensional situation, the diagrams are of a 2-dimensional affine
slice of the fans in question.  The map $\phi$ is the identity in each case.
\begin{figure}[htb]
\includegraphics*{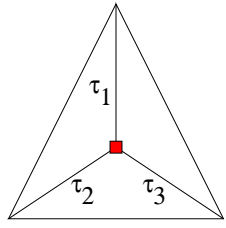}
\caption{The three cases of Lemma 7.1.  The cones 
$\sigma'$ and $\sigma$ are shown in 
grey.}
\end{figure}

\begin{proof}
First, let us 
suppose $\s' \not \in P$.  Then if we are not in (i), there is some
$\tau' \rightarrow \s'$, such that $\tau' \in P$.  Since $\s' \not \in P$
but $\tau' \in P$, $c(\s')=c(\tau')$.  Let us denote $c(\tau')$ by $\tau$.    

The image of $\s'$ is not contained in the boundary of $\tau$, 
since $c(\s')=\tau$.  Thus, we can pick a point $x$ in the interior of
$\s'$ such that $\phi(x)$ is in the interior of $\tau$.  

Now consider a ray in $N_\Q$ from 
a point in the interior of $\tau'$, passing through $x$, and continuing beyond.  Let $y$
be a point just past $x$ on this ray.  $\phi(y) \in \tau$, so by properness,
$y$ is in some cone of $\Sigma'$.  

Now we need a lemma:

\begin{lemma} Let $\phi$ be a proper map from $\Sigma'$ to $\Sigma$.  
Let $\tau$ be a cone of $\Sigma$, of codimension $s$.  Then the inverse image
of $\tau$ is a union of cones of $\Sigma'$ all of codimension $s$.
\end{lemma}
\begin{proof} Suppose some $y$ in the inverse image of $\tau$ is not
contained in any codimension-$s$ cone of $\Sigma'$.  By properness, it
is contained in some cone of $\Sigma'$.  Let $\g'$ be the smallest cone
of $\Sigma'$ containing $y$, and let  the codimension of $\g'$ be $t<s$.

  Pick a small line segment 
$\ell'$ contained in $\g'$ with
midpoint $y$, such that the only point of $\ell'$ contained in $\phi^{-1}(\tau)$
is $y$.  Then $\phi(\ell')$ is also a line segment, which we denote $\ell$.  
The endpoints of $\ell$ are in $c(\gamma')$ but not in $\tau$, and 
$\phi(y) \in \tau$ is a convex combination of the endpoints.  This 
contradicts the assumption that $\Sigma$ is a 
fan.  \end{proof}

From Lemma 7.2 it follows that $y$ must be contained in some cone $\gp$ of 
codimension $k-1$. By construction, $\g'$ has $\s$ in its boundary.   
Also by
Lemma 7.2, $c(\gp)=\tau$, and so $\gp \in P$.    

Now suppose there was some other cone $\delta'\rightarrow \s'$, such that
$\delta' \in P$.  Now $c(\delta') \supset c(\s') = \tau$, so $\delta' \in P$ implies that
$c(\delta')=\tau$.  

Thus $\tau'$, $\g'$, and the putative $\delta'$ are all contained
in the affine span of 
$\phi^{-1}(\tau)$, which is a dimension 
$n'-k+1$ linear subspace of $N'$.  Now $\s'$ is a cone of codimension one in 
that subspace, so it can be contained in at most two larger cones in the
subspace, which means that no cone $\delta'$ can exist.   
Thus, we have established all the conditions of (ii).  

Now suppose $\s' \in P$. 
Let the cones of dimension one higher containing $c(\s')$ be 
$\tau_1,\dots,\tau_k$.  
I claim that for each such $\tau_i$ there is a corresponding
cone $\tau'_i$ of dimension $n'-k+1$ containing $\s'$, such that the image of 
$\tau'_i$ is contained
in $\tau_i$.  

By Lemma 7.2, the inverse image of $\tau_i$ is a union of cones of 
codimension $k-1$.  Pick a point in the interior of $\s'$.  
It is in the inverse image of $\tau_i$, so there   
is some cone $\tau_i'$ of dimension $n'-k+1$ which contains it, and therefore 
all of $\s'$.  

Since $\s'$ is a full-dimensional region in the boundary of the inverse image 
of $\tau_i$, and it fully lies in some $\tau_i'$, it cannot lie in any other 
$(n'-k+1)$-dimensional cone in the inverse image of $\tau_i$.

Thus, the situation is that $\s'$ is contained in cones $\tau_i'$ as above
and also possibly certain other cones, whose images are contained in no
$\tau_i$, and which are therefore not contained in $P$.  Thus we have
established the conditions for (iii).  \end{proof}

Now, we use Lemma 7.1 to prove the theorem.  

In case (i), $f^*(D)\cdot [V(\s')]$ is a linear combination of 
$[V(\tau)]$ with $\tau \rightarrow \s$, but since all such $\tau$ are not
in $P$, $f_*(f^*(D)\cdot [V(\s')])=0$, and since $\s \not \in P$,
$f_*([V(\s)])=0$, so $D \cdot f_*([V(\s)])=0$.  So both sides of the
desired identity are zero.  

Now we consider case (ii).  Again, $D \cdot f_*([V(\s)])=0$.
If the local equation for $D$ on $\tau$ is $m$, the local equation for
$f^*(D)$ on $\tau'$ and $\g'$ is $\phi^*(m)$.  
So the coefficients of
$[V(\tau')]$ and $[V(\g')]$ are $\la \psp(\phi^*(m)),n_{\tau',\s'} \ra$ and 
$\la \psp(\phi^*(m)),n_{\g',\s'} \ra$ respectively.  
  Observe that 
$N_{\tau'}=N_{\g'}$.  Now $n_{\tau',\s'}$ is a generator of $N_{\tau'}/N_{\s'}$,
while $n_{\g',\s'}$ is a generator of $N_{\g'}/N_{\s'}$, but $n_{\tau',\s'}$ 
and $n_{\g',\s'}$ are oriented in opposite directions.  So $n_{\tau',\s'}=-
n_{\g',\s'}$, and therefore $f_*(f^*(D)\cdot[V(\s')])=0$.

Now we consider case (iii).
We wish to show that the coefficients of $[V(\tau_i)]$ in 
$f_*(f^*(D)\cdot [V(\s)])$ and $D \cdot f_*([V(\s')])$ are the same for all
$i$.  For simplicity
of notation, fix an $i$, and let $\tau=\tau_i$, $\tau'=\tau_i'$.  

Let $m$ be the
local equation for $D$ on $\tau$.  
Then $\phi^*(m)$ is 
the local equation for $f^*(D)$ on $\tau'$.  
The coefficient of $[V(\tau')]$ in $f^*(D)\cdot \vz$ is 
$\la \psp(\phi^*(m)),n_{\tau',\s'} \ra$.  We can write (in a unique way)
$m=\ps(m)+u$ for some $u \in \Psi(\s)$.  
Applying $\phi^*$ to this, we get
$\phi^*(m)=\phi^*(\ps(m)) + \phi^*(u)$.  Now $u \in \Psi(\s)$, so
$\phi^*(u) \in \Psi'(\s')$ by the assumed agreement condition between
$\Psi$ and $\Psi'$.  And $\phi^*(\ps(m)) \in M'(\s')$.
So $\psp(\phi^*(m))=\phi^*(\ps(m))$.  So the coefficient of $[V(\tau')]$
in $f^*(D) \cdot \vz$ is $\la \phi^*(\ps(m)),n_{\tau',\s'} \ra =
\la \ps(m),\phi(n_{\tau',\s'}) \ra$.  
When we push this forward, we find
the coefficient of $[V(\tau)]$ in $f_*(f^*(D)\cdot [V(\s)])$ to be
$$[M(\tau'):\phi^*(M(\tau))]\la \ps(m),\phi(n_{\tau',\s'}) \ra.$$

On the other hand, the coefficient of $[V(\tau)]$ in 
$D \cdot f_*([V(\s')])$ is $$[M(\s'):\phi^*(M(\s))] 
\la \ps(m),n_{\tau,\s} \ra.$$

We wish to show that these two expressions are equal. 
Observe that $\phi$ induces a map from $N'(\s')_{\tau'}$ to $N(\s)_\tau$, 
which implies that
$\phi(n_{\tau',\s'})$ is in $N(\s)_\tau$, and thus an integer 
multiple of $n_{\tau,\s}$.  Further, it is a positive integer
multiple, because they
point in the same direction.  

The lattice $M(\s')$ is generated by $M(\tau')$ and one other element, say $u$.
The lattice $\phi^*(M(\s))$ is generated by
$\phi^*(M(\tau))$ and one other element, say $tu$ with $t>0$.  

$\la u,n_{\tau',\s'} \ra =1$, so $\la tu,n_{\tau',\s'} \ra =t$.  But
$\la\phi(tu), n_{\tau,\s} \ra=1$.  So $\phi(n_{\tau',\s'})=tn_{\tau,\s}$.  

But it is also clear that 
$t=[M(\s'):\phi^*(M(\s))]/[M(\tau'):\phi^*(M(\tau))]$.
This establishes the identity in the third case, which completes the proof
of the theorem.  \end{proof}

\section{Ring Structure on $Z_*(X)_\Q$ for Simplicial $\Sigma$}
 
In the next several  sections, we consider the case where $\Sigma$ is 
simplicial, i.e.,~for all $k$, the $k$-dimensional cones of $\Sigma$ 
have $k$ extremal
rays.  Given this condition, the map from $\Q$-Cartier divisors to 
$ Z_{n-1}(X)_\Q$ is an isomorphism.  Thus, we can use $D_i$ to refer both
to the invariant cycle and the associated $\Q$-Cartier divisor.  In this case, as
the title of the section advertises, we get a ring structure on 
$Z_*(X)_\Q=\bigoplus_{k=0}^n Z_k(X)_\Q$.

The definition of the ring structure proceeds in a somewhat unusual order,
as in \cite{DF}.
First, observe that by Theorem 6.1 (Commutativity), 
$\zq$ is a module over $\Q[Y_1,\dots,Y_r]$,
the polynomial ring in 
$r$ variables ($r$ being as always the number of rays of $\Sigma$), where
$Y_i$ acts on $\zq$ by $D_i \cdot{}$.  (The simplicialness of $\Sigma$
is used here to give us that the $D_i$ are $\Q$-Cartier.)

The essential observation is that $\zq$ is actually generated as a module
over $\Q [Y_1,\dots,Y_r]$ by $[X]$.  Let $\tau$ be a cone of dimension
$k$, and renumber the rays of $\Sigma$ so that $\tau$ contains exactly
$\rho_1,\dots,\rho_k$.  Now consider $D_{\rho_1}\cdot\dots\cdot D_{\rho_k} 
\cdot [X]$.  Since the $D_i$ meet properly, this is just their geometric
intersection $[V(\tau)]$, up to multiplication by a non-zero rational.
Thus $[V(\tau)]$ is in the image of $\qyr \cdot [X]$, as desired.  
Since $\zq$ is a cyclic module over $\Q[Y_1,\dots,Y_r]$, $\zq \cong 
\Q[Y_1,\dots,Y_r]/I$ for some ideal $I$.  This puts a ring 
structure on $\zq$.

Naturally, we would like to give explicit generators
for $I$.  
We begin the process in this section.
Let $\alpha$ denote the ring homomorphism $\qyr \rightarrow \zq$.  $I$ is 
the kernel of $\alpha$.  

Define the ideal 
$$I_\Sigma = \la Y_{i_1}\cdots Y_{i_s} | \rho_{i_1},\dots,\rho_{i_s} 
\text{ are not all contained in some cone of } \Sigma \ra.$$
This is a Stanley-Reisner ideal, see \cite{St}.  
For now, we prove the following lemma:

\begin{lemma} $I_\Sigma$ is contained in $I$.  \end{lemma}
\begin{proof} This is a triviality.
Let $Y_{i_1}\cdots Y_{i_s}$ be a generator of $I_\Sigma$.  Its image under 
$\alpha$
 is $D_{i_1} \cdots  D_{i_s}$, which must consist
of a linear combination of $[V(\tau)]$ with $\tau$ containing all the 
$\rho_{i_j}$ --- but there are no such $\tau$ by assumption.  So its image
is zero, and it is in $I$.  \end{proof}

\section {A Presentation for the Ring $\zq$ Defined by an Inner
Product Action}

We can say more about $I$ in the case where $\Psi$ is induced from an inner
product on $M$.

The inner product on $M$ determines an isomorphism $\omega: \Nq \rightarrow
M$ by requiring that
$\la \cdot,\omega(v)  \ra_M
= \la \cdot,v \ra$, where the first set of angle brackets refer to the 
inner product in $M$, and the second to the duality pairing.

Fix $j$, $1\leq j \leq r$.  Let $E_j$ be the principal $\Q$-Cartier
divisor whose local equation on each maximal $\s$ is $\omega(v_j)$.  
I claim that 
$E_j \cdot D_j=0$.    We compute $E_j \cdot D_j$ by viewing $E_j$ as 
a $\Q$-Cartier divisor, and $D_j$ as a cycle.  By the definition of $\omega$, 
$\omega(v_j)$ is perpendicular in $M$ to $v_j^\perp=\rho_j^\perp$.  Thus,
$\pi_{\rho_j}(\omega(v_j))=0$, so the divisor on $D_j$ obtained from
$E_j$ is zero, so $E_j \cdot D_j = 0$.  

We wish to convert this into a statement about $I$, the 
kernel of the map 
$\alpha: \Q [Y_1,\dots,Y_r] \rightarrow \zq$.  To do this, we must write
$E_j$ as a sum of $D_i$.  It is easy to see that  
$$E_j = \sum_{i=1}^r \la\omega(v_j), v_i \ra D_i.$$  

So, if we let $J$ be the ideal in $\qyr$ generated by $\sum_{j=1}^r 
\la \omega(v_j),v_i \ra Y_iY_j$ for $1\leq i \leq r$, then 
$J$ is contained in $I$.  
(The argument above, largely due to Bill Fulton, is a 
simplification of my  original argument.)  

Now we can give a  presentation of $Z_*(X)_\Q$.  

\begin{theorem} When $\Psi$ is induced from an inner product on $M$, 
$I = I_\Sigma + J$, i.e., 
$$  \zq=\qyr/(I_\Sigma + J).$$
\end{theorem}

\begin{proof} As we already remarked, we know that $I_\Sigma + J$ is 
contained in $I$.  To show the other containment, we proceed as in \cite{DF}.  

For $\sigma \in \Sigma$, denote by $Y_\sigma$ the product of those $Y_i$
such that $\rho_i$ is a ray of $\sigma$.  The collection of all the $Y_\s$
form a $\Q$-basis for $\qyr/I$.  So, given any monomial in $\qyr$, it 
suffices to show that it is congruent to a linear combination of $Y_\s$ mod
$I_\Sigma + J$.  

Let $T$ be a monomial in $\qyr$.  Let $|T|$ denote the set of $i$ such that 
$Y_i$ appears in $T$.  Let $|\s|$ denote the set of $i$ such that $\rho_i$
is contained in $\sigma$.  Now, we divide into two cases. First suppose
that $|T|$ does not equal $|\s|$ for any $\s \in \Sigma$.  Then $T$ is
divisible by one of the generators of $I_\Sigma$, so it is congruent to 0
mod $I_\Sigma$, and we are done.  

Now suppose that $|T|=|\s|$ for some $\s$. 
We first prove as a lemma 
the desired reducibility for certain values of $T$, and then use that
to prove the general result.  

\begin{lemma} For $i \in |\sigma|$, $Y_iY_\s$ is congruent to a 
linear combination of $Y_\tau$ mod $I_\Sigma + J$.  \end{lemma}

\begin{proof} Renumbering, we may assume that
the $\rho_j$ which are in $\sigma$ are numbered 1 through $k$.  
Then, for 
$1 \leq j \leq k$, 
$$\sum_{i=1}^{r} \la \om(v_j),v_i \ra Y_iY_\s \in J.$$
Equivalently, for $1 \leq j \leq k$
\begin{equation}
\sum_{1 \leq i \leq k}   \la \om(v_j),v_i \ra Y_iY_\s \equiv 
- \sum_{k+1 \leq i \leq r} \la \om(v_j),v_i \ra Y_iY_\s \mod J. 
\end{equation}
Observe that the monomials on the right-hand side are square-free.

We express (2) in matrix notation.
Let $A$ be the $k \times k$ matrix  whose entries are 
$a_{ji}= \la \om(v_j),v_i \ra$.
Let $W$ be the $k \times 1$ matrix $w_{i1} = Y_iY_\s $.  Let B be the 
$k \times (n-k) $ matrix $b_{ji} = \la \om(v_{j}),v_{i+k} \ra $.
Let $U$ be the $(n-k) \times 1$ matrix $u_{i1}=Y_{i+k}Y_\sigma$.

Then (2) can be rewritten as 
$$ AW \equiv  -BU \mod J. $$

As in \cite{DF}, the critical fact is that the matrix $A$ is invertible.  We prove
this as follows.  Extend the inner product on $M$ to an inner product on
$M_\mathbb{R} = \mathbb{R} \otimes M$.  
Pick an orthonormal basis for $M_\mathbb{R}$ and let $V$
be the $k \times n$ matrix whose $i$-th row is $\omega(v_i)$ written out in 
that
basis.  Then $A = VV^t $.  Since $V$ has 
rank $k$, by the Cauchy-Binet formula $A$ has non-zero determinant, and
is therefore invertible.  

So we can write 
$$ W= -A^{-1}BU \mod J, $$
which gives us an expression for $Y_iY_\sigma$  ($1 \leq i \leq k$)
as a linear combination of 
square-free monomials.  This suffices because every square-free monomial is
either in $I_\Sigma$ or equals $Y_\tau$ for some $\tau \in \Sigma$.
So for $1 \leq i \leq k$, $Y_iY_\s$ is congruent to a linear combination of
$Y_\tau$ mod $I_\Sigma + J$, which proves the lemma.
\end{proof}

Returning to the general case, let $T$ be a monomial in $\qyr$ 
with $|T|=\s$, and let 
$T = \prod_{i\in |\s|} Y_i^{a_i}$.  Define its excess exponent to be 
$e(T)=\sum_{i\in |\s|} (a_i -1)$.  
So a monomial is square-free iff its excess 
exponent is zero.  Now, using Lemma 9.1, it is easy to 
see that we can express $T$ as a sum of terms with smaller excess 
exponents, and this suffices by induction.  
\end{proof}

\section {A Hard Lefschetz-type Theorem for $\zq$ Defined by 
an Inner Product
Action}

{\renewcommand{\omega}{\xi}
Let $\Sigma$ be simplicial.  
In this section we establish a hard Lefschetz-type theorem for $\zq$, where
$\Psi$ is induced from a ``generic'' inner product, an expression
which we 
immediately explain.  By fixing a basis for $M$, we can identify the 
set of inner products on $M$ with the set $P$ of $ n \times n$ symmetric
matrices with coefficients in $\Q$ and all eigenvalues positive.  
The set of $n \times n$ symmetric matrices can be thought of as a finite vector
space over $\Q$, and as such has both a Hausdorff and a Zariski topology. 
$P$ is easily seen to be an open set in the Hausdorff topology.  We say 
that a property holds for a generic inner product if there is a non-empty
Zariski-open subset $S$ of the symmetric matrices such that any matrix in 
$P \cap S$ induces an inner product on $M$ which has the desired property.
(In particular, because of the nature of the Zariski topology, this 
guarantees that there is a non-empty Hausdorff open subset of $P$ such that
any matrix in the subset induces an inner product with the desired property.)

The statement of the theorem is as follows:

\begin{theorem} Let $\Sigma$ be a simplicial fan.    Let 
$\omega = \sum_{i=1}^r a_iD_i$, with all the $a_i$ non-zero.  
Then, for a generic inner 
product on $M$,
 the map from $Z_{n-i}(X)_\Q$ to 
$Z_i(X)_\Q$ induced by multiplication by $\omega^{n-2i}$ is an injection for 
$0 \leq i \leq n/2$.  \end{theorem}

\begin{proof}  Usually a hard Lefschetz theorem would say that $\omega^{n-2i}$
should induce an isomorphism from $Z_{n-i}(X)_\Q$ to $Z_i(X)_\Q$ 
but we cannot 
hope for that because the spaces usually have different dimensions.  If they 
have the same dimensions 
(as for instance in the affine case) then injection and isomorphism are
of course the same thing.  

First, we prove a very special case.  

\begin{lemma} Let $\Sigma$ be affine simplicial, let $\omega = 
\sum_{i=1}^n a_iD_i$, with all the $a_i$ non-zero, and let the inner
product be such that the $v_i$ are all orthogonal to one another.  Then the 
map induced by multiplication by $\omega^{n-2i}$  from $Z_{n-i}(X)_\Q$ to 
$Z_i(X)_\Q$ is an isomorphism. \end{lemma}

\begin{proof}  The choice of inner product means that 
$$\zq \cong \qyn/\la Y_1^2, \dots, Y_n^2 \ra, $$
where the isomorphism takes $D_i$ to $Y_i$.  

Let $W_i = a_iY_i$.  
$\qyn/\la Y_1^2, \dots, Y_n^2 \ra \cong \Q[W_1,\dots,W_n]/\la W_1^2,\dots
, W_n^2 \ra$. 
For $K \subset \{ 1, \dots, n \}$, let $W_K$ be the product of the $W_i$ with
$i \in K$.  Then the $W_K$ form a basis for $\zq$ in which multiplication
is particularly simple: $W_K W_L = W_{K \cup L}$ if $K \cap L = \emptyset$
and 0 otherwise.  In this basis $\omega= W_1 + \dots + W_n$.

To prove the lemma, it suffices to show that for any $i \leq n/2$, and 
any set $K$ of size
$n-i$, $W_K$ is a multiple of $\omega^{n-2i}$.  A quick proof, due to Matt
Frank, is as follows.  We can assume $K = \{ 1,\dots,n-i \}$.  Let
$\zeta = W_{n-i+1}+\dots+W_n$.  Then

$$W_K = \frac{1}{(n-i)!}(\omega - \zeta)^{n-i}= 
\frac{1}{(n-i)!}\sum_{j=0}^{n-i} \binom{n-i}{j} 
\omega^{n-i-j}\zeta^{j}.$$
But for $j>i$, $\zeta^j=0$.  So 
$$W_K= \frac{1}{(n-i)!}\sum_{j=0}^{i}\binom{n-i}{j}\omega^{n-i-j}\zeta^{j}=
\frac{\omega^{n-2i}}{(n-i)!}\sum_{j=0}^i \binom{n-i}j \omega^{i-j}\zeta^j.$$
This proves the lemma. \end{proof}

Keep $\Sigma$ fixed and affine, and consider the set of inner products 
 such that for $0 \leq i \leq n/2$, $\omega^{n-2i}$
induces an isomorphism from $Z_{n-i}(X)_\Q$ to $Z_i(X)_\Q$.  The conditions
which define this subset amount to the non-vanishing of certain determinants,
and are thus 
Zariski-open conditions on the set of symmetric matrices.  
The 
lemma proves that the set is non-empty.  
So there is a non-empty 
Zariski-open subset $S$ of the symmetric matrices such that any element
of $S\cap P$ induces an inner product with the desired property.
Thus, a generic inner product on
$M$ has the property that $\omega^{n-2i}$ induces an isomorphism
between $Z_{n-i}(X)_\Q$ and $Z_i(X)_\Q$.  

Now we return to the proof of the theorem as originally stated, and let 
$\Sigma$ be any simplicial fan.  The condition
that $\omega^{n-2i}$ induces an injection can be checked locally on affine
cones, and thus the set of inner products with the desired property is an 
intersection of non-empty Zariski-open sets, and hence is itself non-empty and
Zariski-open.  \end{proof}
}

\section{An Example: Projective Space}

\begin{con} For the duration of this section,
we rescind our convention that $M$ means $\Qt M$ --- we will do our 
tensoring explicitly.  \end{con}

In this section we consider one way of representing $\mathbb{P}^n$ as a toric
variety $W$, and describe the ring $Z_*(W)_\Q$ which arises from it using
a natural choice of inner product.

Let 
$N'=\mathbb{Z}^{n+1}$ with basis vectors $e_1,\dots,e_{n+1}$.  Let $N =
\{\sum_{i=1}^{n+1} a_ie_i \mid \sum a_i=0\}$.
Let $v_i= e_i-e_{i+1}$, interpreted cyclically, so that $v_{n+1}=e_{n+1}-e_1$. 
Define a fan $\Sigma$ in $N$ consisting
of all cones $\s$ generated by a proper subset of $\{v_1, \dots, v_
{n+1} \}$.  Then $v_i$ is the first lattice point along the ray 
generated by $v_i$, which we call $\rho_i$.
We denote by $D_i$ the Cartier divisor corresponding to $\rho_i$.    
Let $W$ be the toric variety associated to $\Sigma$.  Then 
$W$ is isomorphic to $\mathbb{P}^n$.  

Let $M'$ be dual to $N'$, and $e_i^*$ its dual basis.  Let $M$ be dual to $N$.
Then $M=M'/\mathbb{Z}(e_1^*+\dots+e_n^*)$.  
Let $\bar e_i^*$ denote the image of $e_i^*$
in $M$.  

We put the standard inner product on $M'_\Q$, and identifying $M_\Q$ 
with the orthogonal complement of $e_1^*+\dots+e_{n+1}^*$, 
we induce an inner product on $M_\Q$.  
By the results of the previous section, we see that the ring on the 
invariant cycles is
as follows:
\begin{eqnarray*}
Z_*(W)_\Q &\cong& \Q[Y_1,\dots,Y_{n+1}]/(\la Y_i(Y_i - \frac 12 Y_{i-1} -\frac
12  Y_{i+1}) \mid 1 \leq i \leq n+1 \ra \\
& & \qquad \qquad \qquad \qquad \quad + \la Y_1Y_2\dots Y_{n+1} \ra). \end{eqnarray*}

This can also be seen by direct calculation of $D_i^2$ in $Z_*(W)_\Q$; 
the reader is 
encouraged to try this out, to get the flavour of the computation of intersection
cycles.  

\section{A Cycle-level  Todd Class for Projective Space}

$\Td(X)$, the Todd class of $X$, is an element of $A_*(X)_\Q$.  If $X$ is 
a non-singular toric variety, it satisfies the following formula (see 
\cite{F2}):
$$\Td(X)=\prod_{i=1}^r \frac{D_i}{1-\exp(-D_i)},$$
where the computation takes place in the Chow ring of $X$.  
By interpreting this
as taking place in the ring $Z_*(X)_\Q$, we may obtain a cycle-level 
Todd class for 
any non-singular variety.  For a non-singular variety $X$, let $t(X)$ denote
the cycle-level Todd class obtained in this way.  
Then we have the following theorem. 

\begin{theorem} Let $W$ be the realization of $\mathbb{P}^n$ as a toric variety
discussed
in the previous section.  For $\sigma$ any cone in $\Sigma$, 
the coefficient of $[V(\s)]$ in $t(W)$ is the fraction
of the linear span of $\s$ which is contained in $\s$.  \end{theorem}

\begin{proof} In \cite{DF}, a ring on the cycles of $W$ is constructed 
in which $D_i^2=D_i\cdot D_{i+1}$ (with indices interpreted cyclically), and
it is shown there that the cycle-level  Todd class in this ring has the property 
described in the statement of the theorem.
The essential details of the proof in \cite{DF} go through for our ring without 
any changes. 

Define the formal power series:
\begin{eqnarray*}
\Phi(x)& =& \frac x{1-e^{-x}}\\
Y(x)&= &-\log(1-x)=x+\frac {x^2}2+ \frac {x^3}3 + \cdots .
\end{eqnarray*}

Observe that $ x \Phi(Y(x))=Y(x).$

We consider the set $\{1,\dots,n+1 \}$ arranged around a circle.  
If  $S$ is a subset of 
$\{ 1,\dots,n+1\}$, we call the components of $S$ its maximal 
subsets of consecutive integers (understood cyclically).  
Define $q(S)= \prod \frac 1{|T|+1}$ where
the product is taken over all components $T$ of $S$.  
Let $D_S = \prod_{i \in S} D_i$.

\begin{lemma}
In $Z_*(W)_\Q$, 
$$
t(W)=\prod_{i=1}^{n+1} \Phi(D_i) = \sum_{S \subsetneq \{1,\dots,n+1\}} q(S)D_S.
$$
\end{lemma}

\begin{proof}
We think of the monomials $D_S$ for $S \subsetneq \{1,\dots,n+1\}$ as a 
$\Q$-basis for the ring $Z_*(W)_\Q$.
Let $r(S)$ be the coefficient of $D_S$ in the above product.  Note that since
the monomials produced by multiplying a monomial by $D_i$ all include $D_i$, 
$r(S)$ is
unchanged if we omit all the $\Phi(D_i)$ for $i \not \in S$.  Also, it is
clear that $r(S)=\prod r(T)$ as $T$ ranges over the components of $S$.  
Since this is also true of $q(S)$, it suffices to consider $S$ connected, 
so by symmetry we may assume
that $S=\{1,\dots,k\}$.  By the same argument as above, we may assume that
$n=k$.  

Now consider the map from $Z_*(W)_\Q$ to $\Q[x]$ which takes all the $D_i$ to
$x$.  Then the sum of $q(S)$ over all subsets $S$ of size $n$ of $\{1,\dots,n+1\}$ is 
the coefficient of $x^{n}$ in $\Phi(x)^{n+1}$.  

Given a power series $f(x)$, denote the $x^n$ coefficient of $f(x)$ by 
$[x^n]f(x)$.
Then by the Lagrange Inversion Formula (see \cite{Wi}),
$$ (n+1) [x^{n+1}] Y(x)= [x^n]\Phi(x)^{n+1}.$$

So $[x^n]\Phi(x)^{n+1}=1$.  So the sum over all the subsets $S$ of size $n$
of $\{1,\dots,n+1\}$ is 1.  By symmetry, then, for any such $S$, 
$r(S)=\frac 1{n+1}$, so $r(S)=q(S)$ as desired.  \end{proof}

Thus, for a cone $\s \in \Sigma$, define $S$ to be the indices of the rays
of $\Sigma$ contained in $\s$.   Then
the coefficient of $[V(\s)]$ in $t(W)$ is $q(S)$, and, as stated in \cite{DF},
this is the fraction of the linear span of $\s$ which is contained in
$\s$, which proves the theorem.  \end{proof}

\section{Computing the Coefficient of $[V(\s)]$ in $D^n$}

In this section, we compute the coefficient of $[V(\s)]$ in $D^n$, where
$\s$ is an $n$-dimensional cone of $\Sigma$, and express it as the sum
of plus or minus the volume of a large number of simplices.  The result 
is somewhat unwieldy, but gives the basic geometrical picture which we 
shall exploit in the following two sections.  

We begin with some very elementary statements about volumes.  
Given a simplex $S$ with vertices $u_0, \dots, u_n$ in $M$, 
consider the matrix $A$ whose $i,j$-th entry is given by
$\la u_j - u_0,v_i \ra$, where the $\{v_j\}$ are a basis of $N$. 
Then the lattice volume of 
$S$, written $\avol(S)$ or $\avol(u_0,u_1,\dots,u_n)$
is $\frac{1}{n!}|\det(A)|$.  This is a well-defined element of $\Q$, and
doesn't depend on the order of the vertices.  In general, when we refer to
volume, we mean volume normalized in this way with respect to the lattice.
Note that this is only possible if the affine span of the set whose volume
we wish to take is rational, i.e.~is parallel to some sublattice of
the group $M$ (not meaning $\Qt M$ here), 
so that there is a well-defined lattice
in the affine span of the set.    

\begin{lemma}
Let $\rho$ be a rational ray in $N$, $v$ the first lattice point along it.
Suppose $u_0, \dots, u_{n-1} \in \rho^\perp$.  Then 
$$\avol(u_0,u_1,\dots,u_n) =
\frac 1n \la u_n,v \ra \avol(u_0,\dots,u_{n-1}).$$   \end{lemma}

\begin{cor*} Let $\rho$ be a rational ray in $N$, $v$ the first
lattice point along it.  Let $P$ be a pyramid in $M$ over a polytope $Q$ in $\rho^\perp$,
with apex $u_n$.  Then 
$$ \avol(P) = \frac 1n \la u_n,v \ra \avol(Q).$$ \end{cor*}

Let $D$ be a $\Q$-Cartier divisor on $X$, a general $n$-dimensional
toric variety.  Let $\s$ be an $n$-dimensional cone in $\Sigma$, the fan
corresponding to $X$.  We will compute the coefficient of $[V(\s)]$ in 
$D^n$.  Let the local equation for $D$ on $\s$ be $m \in M$.  

Clearly, the coefficient of $[V(\s)]$ in $D^n$ is

$$ \sum_{0=\g_0\prec\g_1\prec\dots\prec\g_n=\s} \,\, \product_{1 \leq i\leq n}
\la \p_{\g_{i-1}}(m),n_{\g_i,\g_{i-1}} \ra $$

Now we analyze the contribution from each term in the sum.  

\begin{lemma} Fix a chain $0=\g_0\prec\g_1\prec\dots\prec\g_n=\s$. 
The absolute value of the contribution from this chain to the coefficient of 
$[V(\s)]$, 
$$\left| \product_{1 \leq i\leq n} \la \p_{\g_{i-1}}(m),n_{\g_i,\g_{i-1}} 
\ra \right|$$ 
is $n!$ times 
the lattice volume
of the simplex with vertices $\pi_{\g_0}(m)$, $\pi_{\g_1}(m)$,$\dots$,
$\pi_{\g_n}(m)$.  \end{lemma}

\begin{proof}
The proof is by induction on $n$.  For $n=1$ it is clear.  Suppose it is
true for $n-1$.  Then $\left| \product_{1 \leq i\leq n} 
\la \p_{\g_{i-1}}(m),n_{\g_i,\g_{i-1}} \ra \right|$ is $(n-1)!$ times the
lattice volume of the simplex with vertices 
$\pi_{\g_1}(m),\dots,\pi_{\g_n}(m)$.  

We apply Lemma 13.1, with $\g_1$ for $\rho$. Observe that $v$ is $n_{\g_1,
\g_0}$.  The result follows.  \end{proof}

Suppose that $X$ is a proper toric variety, and $D$ a 
$\Q$-Cartier divisor on it
with local equations $m_\s$.  Let $S$ be the set of all
$-m_\s \in M$, for $\s$ a maximal cone.  Form the convex hull of $S$, and
call it $P$.  If there is a 1-1 inclusion-reversing
map $F$ from the cones of $\Sigma$ to
the faces of $P$, which extends the map from the maximal cones to $S$, we say
that the divisor $D$ corresponds to the polytope $P$.

In this case, we see that if we add together the simplices of Lemma 13.2, taking
care to cancel overlaps with opposite signs, we obtain $P$, from which we may
recover the well-known fact that
for a Cartier divisor
corresponding to $P$, $\deg(D^n)=n!\avol(P)$.

\section{Computing $D^n$ 
for the Generic Flag Action}

In this section, we investigate the coefficient
of $[V(\s)]$ in $D^n$ for $\Psi$ given by a generic flag in $M$, 
under the assumption
that $\s$ is simplicial.  

Let $0=L_0 \subset L_1\subset \dots \subset L_n=M$ be a 
complete generic flag
in $M$, and $\Psi$ the corresponding choice of complements.  
Let $\s$ be an $n$-dimensional simplicial cone.  
Let $D$ be a $\Q$-Cartier divisor,
and let its local equation on $\s$ be $m$.  Let the facets of $\s$ be
$\tau_1,\dots,\tau_n$, and the extreme rays of $\s$ be $\rho_1,\dots,\rho_n$,
with $\rho_i$ the unique extreme ray not contained in $\tau_i$.    
Let $\sigma^\vee$ be the dual cone to $\sigma$ --- by definition this
is the cone consisting of elements of $M$ whose duality pairing with all the
elements of $\sigma$ is non-negative.   
For
$\g \prec \s$, let 
$\g^* = \g^\perp \cap \sigma^\vee$.  
The $\tau_i^*$
are the extreme rays of the cone $\sigma^\vee$, and the $\rho_i^*$ are 
the facets of $\s^\vee$.  

A simplex in $M$ determines $n+1$ 
hyperplanes, the affine spans of its facets.  
Conversely, hyperplanes $H_1,\dots,H_{n+1}$ 
in general position determine a simplex, as follows.  
Let $v_i$ be the intersection of all the hyperplanes
excluding $H_i$.  Then we say that the collection of hyperplanes determines
the simplex which is the convex hull of the $v_i$.  

We prove the following 
theorem:

\begin{theorem} Let $\s$ be an $n$-dimensional simplicial cone.
The coefficient of $[V(\s)]$ in $D^n$ is plus or minus
$n!$ times the volume of the simplex corresponding to the affine spans of the 
$\rho_i^*$ and the hyperplane $H$ passing through $m$ parallel to $L_{n-1}$.  
It is positive or negative as the number of $\tau_i^*$ which $H$ does not 
intersect
is even or odd.  \end{theorem}

\begin{proof} The proof is by induction on $k$.  The induction claim is that
if $\delta$ is an $n-k$-dimensional face of $\s$ (so, by renumbering, we may
assume that $\delta$ is the intersection of $\tau_1, \dots, \tau_k$), then

$$\frac 1{k!}
\sum_{\delta=\g_{n-k}  \leftarrow \dots \leftarrow \g_{n}=\s}
\,\, \product_{n-k+1 \leq i\leq n}
\la \p_{\g_{i-1}}(m),n_{\g_i,\g_{i-1}} \ra$$
is plus or minus the volume of the simplex in $\delta^\perp$ 
corresponding to the hyperplanes $\delta^\perp \cap \rho_1^*, \dots
, \delta^\perp \cap \rho_k^*, \delta^\perp \cap H$, signed positive or 
negative as
the number of $\tau_i^*$ for $1 \leq i \leq k$ that $\delta^\perp \cap H$ 
misses is even or odd. 

The 
case $k=0$ is true by convention, and $k=n$ is the desired result.  
Assume 
the statement holds 
for $k-1$.  Then

{\setlength{\arraycolsep}{1mm}
\begin{eqnarray}
\lefteqn{\frac 1{k!}
\sum_{\delta=\g_{n-k} \leftarrow  \dots \leftarrow \g_{n}=\s}
\,\, \product_{n-k+1 \leq i\leq n}
\la \p_{\g_{i-1}}(m),n_{\g_i,\g_{i-1}} \ra }\nonumber \\
& &\!\!\!\!\! = \!\frac 1k \sum_{i=1}^k 
 \left( \!\!\frac{\la \pi_\delta(m), n_{\delta+\rho_i,\delta} \ra}
{(k-1)!}\sum_{\substack{\delta+\rho_i=\g_{n-k+1} \leftarrow \\
\dots \leftarrow \g_n=\s}} \,\,
\product_{j=n-k+2}^n
\!\!\!\!\!\!\la \p_{\g_{j-1}}(m),n_{\g_j,\g_{j-1}} \ra\!\! \right) 
\end{eqnarray}}

Let $S$ be the simplex in $\delta^\perp$ corresponding to the hyperplanes
$\delta^\perp \cap \rho_1^*,\dots,\delta^\perp \cap \rho_k^*$ and 
$\delta^\perp \cap H$.
Let $F_i$ be the face of $S$  corresponding to the hyperplane
$\delta^\perp \cap \rho_i^*$.  Then by the induction
hypothesis, for each $i$, $1\leq i\leq k$,

$$\frac1{(k-1)!}\left| \sum_{\substack{\delta+\rho_i=\g_{n-k+1} \leftarrow \\ \quad \dots \leftarrow \g_{n}=\s}}
\,\,\product_{n-k+2 \leq i \leq n}
\la \p_{\g_{i-1}}(m),n_{\g_i,\g_{i-1}} \ra \right|= \avol(F_i).$$

I claim that $\pi_\delta(m) \in H \cap \delta^\perp$.  By construction it is 
in $\delta^\perp$.  To show it is in $H$, first observe that by the definition
of $H$, $m \in H$.  But now
by construction, $\pi_\delta(m)-m \in \Psi(\delta)=L_{n-k}$.  
$L_{n-k} \subset L_{n-1}$ except when $k=0$, but then $\pi_\delta(m)-m=0$.  
So $\pi_\delta(m)-m \in L_{n-1}$.  So $\pi_\delta(m) \in H$, and hence in 
$H \cap \delta^\perp$.  

Let $C_i$ be the pyramid over $F_i$ with apex $\pi_\delta(m)$.  Then

$$\avol(C_i) = \frac 1k | \la \pi_\delta(m),n_{\delta+\rho_i,\delta} \ra|\avol(F_i).$$

So (3) is a sum of plus or minus
the volumes of the $C_i$.  Let $s \in \{1,-1\}$ be positive if 
the number of $\tau_i^*$ with $1\leq i \leq k$ missed by $H$ is even, and
negative if the number missed is odd.  Let $s_i$ denote the sign 
such that
$$\frac 1{k!}
\la \pi_\delta(m),n_{\delta+\rho_i,\delta} \ra
\sum_{\substack{\delta+\rho_i=\g_{n-k+1} \leftarrow \\ \quad \dots \leftarrow \g_{n}=\s}}
\,\,\product_{n-k+1 \leq i\leq n}
\!\!\!\!\!\la \p_{\g_{i-1}}(m),n_{\g_i,\g_{i-1}} \ra= s_i\avol(C_i).$$  

Observe that
$\la \pi_\delta(m),n_{\delta+\rho_i,\delta}\ra$ is positive if 
$\pi_\delta(m)$ is on the positive side of $\delta^\perp \cap \rho_i^*$, 
that is,
the same side as $\tau_i^*$.  So $ss_i$ is positive if $H$ hits $\tau^*_i$
and $\pi_\delta(m)$ is on the $\tau_i^*$ side of $\delta^\perp \cap \rho_i^*$
or if $H$ misses $\tau_i^*$ and $\pi_\delta(m)$ is on the 
$-\tau_i^*$ side of $ \delta^\perp \cap \rho_i^*$.  But this amounts to saying 
that $ss_i$ is positive iff $\pi_\delta(m)$ is on the same side of $\delta^\perp
\cap \rho_i^*$ as the intersection of $H$ with the linear span of $\tau_i^*$.  

Now there are two cases.  The first case is if $\pi_d(m)$ is in $S$.  Then
for all $i$,  $ss_i=1$.  Then the cones $C_i$ form a dissection of $S$,
that is to say, their union is $S$, and they overlap only on boundaries.  
So 
\begin{eqnarray*}
\frac 1{k!}
\sum_{\substack{\delta=\g_{n-k} \leftarrow \g_{n-k+1} \leftarrow 
\\ \quad \dots \leftarrow \g_{n}=\s}}
\,\,\product_{n-k+1 \leq i\leq n}
\!\!\!\!\!\la \p_{\g_{i-1}}(m),n_{\g_i,\g_{i-1}} \ra &=& 
\sum_{i=1}^k s_i \avol(C_i) \\
& =& s\avol(S). \end{eqnarray*}

The second case
 is if $\pi_\delta(m)$ is not in $S$.  Then let $Q$ be the convex hull
of $S$ and $\pi_\delta(m)$.  We see that the $C_i$ with $ss_i$ positive form 
a dissection of $Q$, while the $C_i$ with $ss_i$ negative form a 
dissection of $Q \setminus S$.  So
\begin{eqnarray*}
& &\frac 1{n!}
\sum_{\delta=\g_{n-k} \leftarrow \g_{n-k+1} \leftarrow \dots \leftarrow \g_{n}=\s}
\,\,\product_{n-k+1 \leq i\leq n}
\!\!\!\!\!\!\!\la \p_{\g_{i-1}}(m),n_{\g_i,\g_{i-1}} \ra \\ 
& & \qquad \qquad =\sum_{i=1}^k s_i \avol(C_i) \\
& & \qquad \qquad =s(\avol(Q)-\avol(Q \setminus S)) \\
& &\qquad \qquad =s(\avol(S)). \end{eqnarray*} 

This finishes the induction step, and the theorem.  
\end{proof}

\section {General Computations for the Generic Flag Action}

In this section we make use of the result of the previous section to obtain
expressions for products of any number of $\Q$-Cartier divisors in 
$Z_*(X)_\Q$ with $\Psi$ given by a generic flag
 $L_0 \subset L_1 \subset \dots \subset L_n=M$.

First, we reformulate the theorem of the previous section algebraically.

\begin{theorem}
Let $\s$ be an $n$-dimensional simplicial 
cone in $\Sigma$, $\tau_1,\dots,\tau_n$ the facets of $\s$, 
and $D$ a $\Q$-Cartier 
divisor, with local equation $m$ on $\s$. 
Let $w \in N_\Q$ be a non-zero vector 
perpendicular to $L_{n-1}$. 
 For 
$1\leq i \leq n$, pick $t_i \in \tau^*_i$, such that the 
lattice volume of the 
simplex with vertices $0,t_1,\dots,t_n$ is $\frac 1{n!}$.  

Then the coefficient of 
$[V(\s)]$ in $D^n$ is

$$\frac {\la m,w \ra^n}{\la t_1,w \ra \cdots \la t_n,w \ra}.$$
\end{theorem}

\begin{proof}
Consider how much we would need to stretch the simplex 
with vertices located at $0,\:t_1,\dots,t_n$ along 
each of the $\tau^*_i$ in order for it to coincide with the simplex $S$ from
the previous section.  The scaling factor along 
$\tau^*_i$ is $\la m,w \ra/\la t_i,w\ra$.
Thus, the volume of $S$ is:

$$ \avol(S)=\frac 1{n!}
\left| \frac {\la m,w \ra^n}{\la t_1,w \ra \cdots \la t_n,w \ra}
\right|. $$

Now, observe that the sign on the expression inside the absolute value signs
is exactly that dictated by Theorem 14.1, which proves the theorem.  
\end{proof}

From this, we derive a corollary which holds for products of different
Cartier divisors:

\begin{cor} Let $\s$ be an $n$-dimensional 
simplicial cone in $\Sigma$
and let $E_1,\dots$, $E_n$ be $\Q$-Cartier divisors, with local equations
$m_1,\dots,m_n$ on $\s$.  Pick $w$ and $t_1$,$\dots$, $t_n$ as in the 
theorem. Then the coefficient of
$[V(\s)]$ in $E_1 \cdot \cdots  E_n$ is 

$$ \frac {\la m_1,w \ra\cdots \la m_n,w \ra}{\la t_1,w \ra \cdots \la t_n,w \ra}.$$
\end{cor}

\begin{proof}
First, by a straightforward algebraic identity, 
$$ E_1 \cdot \cdots  E_n = \frac 1{n!} \sum_{j=0}^n \sum_{\substack{ S 
\subset \{1, \dots, 
n \} \\ |S|=j}} (-1)^{n-j}(\sum_{i\in S} E_i)^n.$$

Thus, the coefficient of $[V(\s)]$ in $E_1 \cdot \cdots  E_n$ is
\begin{eqnarray*}
& \frac {\displaystyle 1}{\displaystyle\la t_1,w \ra \cdots \la t_n,w \ra}
\sum_{j=0}^n  \sum_{\substack{ S    
\subset \{1, \dots, 
n \} \\ |S|=j}} (-1)^{n-j} \left(\sum _{i \in S}\la m_i,w \ra\right) ^n
 \\
&\qquad=\frac{\displaystyle\la m_1,w\ra \cdots \la m_n,w \ra}
{\displaystyle\la t_1,w\ra \cdots \la t_n,w \ra}. \end{eqnarray*}
This proves the corollary. \end{proof}

Now suppose that we are interested in some degree $k$ polynomial in a 
collection of Cartier divisors $E_1,\dots,E_s$.  We have the following 
corollary:

\begin{cor} Let $\s$ be a $k$-dimensional simplicial cone in
$\Sigma$, let $E_1,\dots,E_s$ be $\Q$-Cartier divisors with local
equations $m_i$ on $\sigma$, and let
$q_k(x_1,\dots,x_s)$ be a homogeneous
polynomial of degree $k$.
Let $w$ be a non-zero vector in $N_\s$ in the direction perpendicular to
$L_{k-1}$. 
Let $t_1, \dots,t_k$ be chosen                      in $M_\s$ to lie on the
extreme rays of the image of $\sigma$ in $M_\s$, so that the lattice volume
of the simplex with vertices $0,t_1,\dots,t_k$ is $1/k!$.

Then the coefficient of $[V(\s)]$ in $q_k(E_1,\dots,E_s)$ is

$$\frac {q_k(\la m_1,w \ra,\dots,\la m_s,w \ra)}
{\la t_1,w \ra \cdots \la t_k,w \ra}.$$
\end{cor}

\begin{proof} The proof is an application of the technique outlined in section
5 to reduce the computation of the coefficient of an arbitrary cone to
that of a cone of full dimension.  \end{proof}

The results in this section and the previous one are only directly 
applicable to computing the coefficient of $[V(\s)]$ 
in a product of $\Q$-Cartier divisors when $\s$ is simplicial.  
However, they could still
be applied to computing the coefficient of $[V(\s)]$ with
$\s$ non-simplicial, as follows.

Suppose $X$ is a non-simplicial toric variety.  Find a 
simplicial toric variety $X'$ and a proper birational map 
$f: X' \rightarrow X$.  (This amounts to choosing a simplicial subdivision
of $\Sigma$.)
Let $E_1,\dots,E_k$ be $\Q$-Cartier divisors on
$X$.  Then by Theorem 7.1,
$$E_1 \cdot \cdots  E_k \cdot [X] = f_*(f^*(E_1) \cdots  f^*(E_k)
\cdot [X']),$$
and 
$f^*(E_1)\cdot \cdots  f^*(E_k)
\cdot [X']$ can be computed using the techniques of this section and the
previous one.  

\section{Canonical Constructions}

\begin{con} For this section we revoke the convention of tacitly 
tensoring $M$ and its sub- and quotient lattices by $\Q$.  \end{con}

 From either the inner product action or the 
generic flag action there is a way to obtain an action of Cartier divisors
on cycles which does not depend on any choice, but one has to allow the
cycles to take coefficients in a much larger field.  The construction we
give is based on a construction in \cite{Mo}.

First, we discuss the canonical inner product action.  Let $B(M)$ be the 
set of symmetric bilinear (not necessarily non-degenerate) forms on $M_\Q$.
$B(M)$ can be thought of as a variety defined over the rationals.  Let 
$\K$ be its rational function field.  Let $B^+(M)$ denote the positive definite
forms in $B(M)$.  
Now, we wish to define an canonical action of 
Cartier divisors on cycles with coefficients in $\K$.

For $D \in \Div(X)$ and $\sigma$ an $n-k$ dimensional cone in $\Sigma$,
we want to define $D \cdot_{\text{IP}} [V(\s)]$ in $Z_{k-1}(X)_\K$, which we
can think of as a rational function from $B(M)$ to $Z_{k-1}(X)_\Q$.  Since
$B^+(M)$ is Zariski dense in $B(M)$, to define a rational function on
$B(M)$ it suffices to define it on $B^+(M)$.  So, for $b\in B^+(M)$, let
$(D\cdot_{\text{IP}} [V(\s)])(b)=D\cdot_b [V(\s)]$, where $\cdot_b$ denotes
the inner product action relative to the inner product given by $b$.  One
easily sees that this is a rational function, and thus we can let this 
consitute the 
definition of $D\cdot_{\text{IP}}[V(\s)]$.  
Extending linearly, we obtain an action of $\Div(X)$ on $Z_*(X)_\K$.

This action has all the properties we proved for the inner product action
(commutativity, agreement with usual intersection modulo rational equivalence,
etc.) because all these properties hold for all $b\in B^+(M)$, which
is Zariski dense in $B(M)$.  

An interesting result of this construction is that it allows us to define
an action of Cartier divisors on cycles with respect to a symmetric 
bilinear form
which is not positive definite: simply evaluate the canonical action at 
that bilinear form.  Note that this will not necessarily be well-defined,
however, since the coefficients in question may blow up there.

The canonical flag action is defined in much the same way.  The set $\Fl(M)$ 
of complete flags in $M_\Q$ is a variety, and we let $\F$ be its function 
field.  For $\Cal{F} \in \Fl(M)$, we define
$$(D\cdot_{\text{Fl}} [V(\s)])(\Cal{F}) =D\cdot_\Cal{F} [V(\s)].$$

This will only be well-defined for $\Cal{F}$ generic with respect to $\sigma$,
but this is an open set in $\Fl(M)$, and therefore dense, so the same
arguments go through.  Thus, we have defined an action of Cartier 
divisors on cycles with coefficients in $\F$, and, as above, it has all the
properties we are used to.  

As an example, we state the result about the canonical flag action
corresponding to Corollary 2 to Theorem 15.1.  

\begin{theorem} 
Let $\s$ be a $k$-dimensional cone.
$E_1,\dots, E_s$ Cartier divisors on $X$, and $q_k(x_1,\dots,x_s)$ a 
homogeneous polynomial of degree $k$. 
Let the local equation for $E_i$ on
$\tau$ be $m_i \in M_\s$.    
Let $t_1,\dots,t_k$ be as the 
statement of Corollary 2.      
Let $w: \Fl(M) \rightarrow \Proj(\Qt N_\s)$, where if $\Cal{F}=F_0\leq \dots
\leq F_n$, then 
$w(\Cal{F})$ is the direction
in $\Qt N_\s$ perpendicular to $F_{k-1}$.

Then the coefficient of $[V(\s)]$ in $q_k(E_1,\dots,E_s)$ is 
\begin{equation*} 
\frac{q_k(\la m_1,w \ra, \cdots ,\la m_s,w \ra)}
{\la t_1,w \ra \cdots \la t_k,w \ra}. \end{equation*} \end{theorem}

\section{Characteristic Classes of Toric Varieties}

Let $X$ be a simplicial toric variety.  Suppose we have a characteristic
class of $X$ which is written as a polynomial $p(D_1,\dots,D_r)$ interpreted 
in
$A_*(X)$.  
(For example,
the $j$-th Chern class of $X$ can be written as the $j$-th symmetric 
polynomial on the $D_i$, and one can use this to obtain a polynomial in 
the $D_i$ representing any class written as a 
polynomial in the Chern classes.) If we pick $\Psi$ and 
interpret $p(D_1,\dots,D_r)$ in 
$Z_*(X)_\Q$, we determine a cycle which represents the characteristic class
in question.
If $\Psi$ was chosen using the generic flag action or the canonical flag
action, we can compute this cycle quite explicitly using 
Corollary 2 to Theorem 15.1 or Theorem 16.1, respectively.
We recover in this case the cycle-level characteristic classes
obtained by Morelli in \cite{Mo} using the Baum-Bott residue formula,
but our result is somewhat more general: the results of \cite{Mo} apply
only when $X$ is non-singular and compact.  

If $X$ is non-singular, one can use this approach to compute a cycle-level
Todd class, as we have already seen in Section 12.  
(Determining cycle-level Todd classes is a matter of some interest, because
of connections to counting lattice points in polytopes, see \cite{Da,F2,Mo}.) 

For singular $X$, $\Td(X)$ is computed by taking a resolution of
singularities $X' \rightarrow X$, computing the Todd class of $X'$, and 
pushing it forward to $X$.  This procedure can be done for the cycle-level 
Todd class of a toric variety as well, 
but {\em a priori} this introduces a dependance on 
the choice of resolution.  In fact, under reasonable conditions,
the resulting cycle-level Todd class does not depend on the choice of resolution,
but the proof of this relies on results from \cite{Po} and 
is somewhat complicated, see \cite{T}.  

There is more work to be done
when it comes to understanding cycle-level Todd classes of toric varieties.
One would like to know, for instance, for what toric varieties Theorem 12.1
holds.  Examples in \cite{Mo} show         that it does not hold for all
toric varieties.  On the other hand
there is reason to think that it may hold for 
toric varieties arising from partially ordered sets, see \cite{T2}.  

\section*{Acknowledgements}

This paper is based on my thesis \cite{T}, written at the University of Chicago
under the direction of Bill Fulton.  I would like to thank him for his help
in shaping its substance and its form.  I would also like to thank Matt
Frank and Ian Robertson for useful conversations.

\end{document}